\documentclass[11pt,leqno]{amsart}
\usepackage{amsmath,stmaryrd}
\usepackage{amssymb}
\usepackage{amsthm}
\usepackage{epsf}
\usepackage{graphicx}
\usepackage{esint} 
\usepackage{dsfont}
\usepackage[pagebackref]{hyperref} 
\hypersetup{pdfpagemode=FullScreen,  colorlinks=true, citecolor=blue} 
\usepackage{enumerate} 

\numberwithin{equation}{section}

\DeclareMathOperator{\diver}{div}

\DeclareMathOperator{\dist}{dist}
\DeclareMathOperator{\supp}{supp}
\DeclareMathOperator{\diam}{diam}

\newcommand{\ms}{\medskip}

\newcommand{\R}{\mathbb R}

\newcommand{\bp}{\noindent {\em Proof: }}
\newcommand{\ep}{\hfill $\square$ \medskip}

\newcommand{\E}{\mathcal E}

\theoremstyle{plain}
\newtheorem{theorem}[equation]{Theorem}

\newtheorem{corollary}[equation]{Corollary}
\newtheorem{proposition}[equation]{Proposition}

\newtheorem{definition}[equation]{Definition}

\theoremstyle{definition}

\theoremstyle{remark}

\begin{document}

\title{Carleson perturbations of locally Lipschitz elliptic operators.}

\author[Feneuil]{Joseph Feneuil}
\address{Universit\'e Paris-Saclay, Laboratoire de math\'ematiques d’Orsay, 91405, Orsay, France}
\email{joseph.feneuil@universite-paris-saclay.fr}

\maketitle

\begin{abstract} 
In one-sided Chord-Arc Domains $\Omega$, we demonstrate that the $A_\infty$-absolute continuity of the elliptic measure with respect to the surface measure remains stable under $L^2$ Carleson perturbations. This stability holds provided that either the elliptic operator $L_0=-\diver A_0\nabla$, which is being perturbed, or the perturbed operator $L_1=-\diver A_1\nabla$ satisfies the condition $\sup_{X\in \Omega }\dist(X,\partial \Omega)|\nabla A_i(X)| <\infty$ on its coefficients. $L^2$ Carleson perturbations are slightly more general than those previously discussed in the literature.

The proof hinges on the availability of a comprehensive elliptic theory and a domain $\Omega$ that allows uniform non-tangential access to any point on its boundary. Consequently, while the current theory of $L^2$ Carleson perturbations can be extended to more general contexts, we have chosen not to do so in order to simplify the presentation.
\end{abstract}

\ms\noindent{\bf Keywords:} Boundary value problems, Elliptic operators with rough coefficients, Carleson perturbations, Reverse H\"older inequality.

\ms\noindent

\tableofcontents

\section{Introduction}
\label{S1}

The Dirichlet problem for elliptic operators with data in \(L^p\) has been extensively studied over the decades and is now well understood. Building on numerous previous works, and assuming {\em a priori} the \((n-1)\)-Ahlfors regularity of the boundary \(\partial \Omega\) and the existence of Corkscrew points, the authors of \cite{AHMMT} provided a complete characterization of the domain \(\Omega \subset \mathbb{R}^n\) for which the Dirichlet problem for the Laplacian is solvable in \(L^p\) for some \(p \in (1, \infty)\). Another interesting and related aspect is the role of the elliptic operator in this characterization, which raises two natural questions: How is the solvability of the Dirichlet problem affected by perturbations? And what are the `optimal' perturbations allowed in the operators? The article \cite{FKP91} addresses this question when the domain is a ball and has since been generalized to weaker settings: Chord Arc Domains (\cite{MPT}), $1$-sided Chord Arc Domains (\cite{CHM, CHMT}), to domains satisfying a capacity condition and Carleson conditions that depend on operators (\cite{AHMT}, higher and mixed codimensions \cite{MP, FP}, and operators with BMO antisymmetric coefficients (\cite{DSU}).

We present here an alternative notion of Carleson perturbations, which is slightly more general than those in \cite{FKP91} and subsequent articles but requires one of the two operators - the one being perturbed or the perturbed one - to satisfy extra regularity on its coefficients. This alternative result has been used in \cite{FenReg} and we were asked to provide details of the proof.

We will discuss the literature and our results in domains known as 1-sided CAD. All theorems presented here, both existing and new, hold in the setting introduced in \cite{DFMhighcod}, and part of them hold in the general context of \cite{DFMmixed} or in the abstract setting of \cite{FP}. However, we have chosen to focus on 1-sided CAD to keep our article moderately concise and to use conditions that the reader is hopefully familiar with. The proofs in a more general context will be included in the book \cite{FMbook} in preparation. Alternatively, the reader may simply prefer to adapt the proofs of \cite{FP} and \cite{FMZ}, as our proof follows the ideas presented in \cite{FP,FMZ}.

\subsection{The $A_\infty$-absolute continuity of the elliptic measure} 
We postpone the presentation of our geometric setting, which is not central to our paper, to Subsection \ref{SsCAD}. Let us introduce the PDE tools necessary for our results.

\begin{itemize}
\item An operator \( L = -\diver A \nabla \) is called \textbf{uniformly elliptic} if there exists a constant \( C > 0 \) such that
\begin{equation} \label{defellip} 
|A(X)\xi \cdot \zeta| \leq C |\xi| |\zeta| \quad \text{and} \quad A(X)\xi \cdot \xi \geq C^{-1}|\xi|^2 \quad \text{for all } X \in \Omega, \, \xi, \zeta \in \mathbb{R}^n.
\end{equation}
A \textbf{weak solution} to \( Lu = h \) in \( D \subset \Omega \) is a function in \( W^{1,2}_{loc}(D) \) that satisfies the weak condition
\begin{equation} \label{defweaksol}
\iint_D A \nabla u \cdot \nabla \varphi \, dX = \iint_{\Omega} h \varphi \, dX \quad \text{for all } \varphi \in C^\infty_0(D).
\end{equation}

\item The \textbf{Green function} \( G_L^Y \) associated with a uniformly elliptic operator \( L \) is a collection \(\{G_L^Y\}_{Y\in \Omega}\) of functions on \(\Omega\) such that for any \( h \in C_0^\infty(\Omega) \), the function defined by
\begin{equation} \label{defGLY}
u^h(X) := \iint_{\Omega} G_L^Y(X) h(Y) \, dY
\end{equation}
is the bounded weak solution of \( Lu = h \) in \(\Omega\) that lies in \( C^0(\overline{\Omega}) \) and satisfies \( u_g|_{\partial \Omega} = 0 \).

\item The \textbf{elliptic measure} \( \omega_L^X \) associated with a uniformly elliptic operator \( L \) is the collection of probability measures \(\{ \omega^X_L \}_{X \in \Omega}\) on \(\partial \Omega\) such that, for any \( g \in C^0_0(\partial \Omega) \), the function defined by
\begin{equation} \label{defhm}
u_g(X) := \int_{\partial \Omega} g(y) \, d\omega^X_L(y)
\end{equation}
is the bounded weak solution to \( Lu = 0 \) in \(\Omega\) that lies in \( C^0(\overline{\Omega}) \) and satisfies \( u_g|_{\partial \Omega} = g \).

\item If \(\mu, \nu\) are two measures on \(\partial \Omega\), we say that \(\mu \in A_\infty(\nu)\) - or that \(\mu\) is \(\mathbf{A_\infty}\)-\textbf{absolutely continuous} with respect to \(\nu\) - if for any \(\epsilon > 0\), there exists a \(\eta > 0\) such that for any ball \( B \) centered on \(\partial \Omega\), and any Borel set \( E \subset \Delta := B \cap \partial \Omega \),
\begin{equation} \label{defAinfty}
\frac{\mu(E)}{\mu(\Delta)} < \eta \implies \frac{\nu(E)}{\nu(\Delta)} < \epsilon.
\end{equation}
It is important to note that the \( A_\infty \)-absolute continuity is an equivalence relationship, that is,
\[ \mu \in A_\infty(\nu) \iff \nu \in A_\infty(\mu) \]
and
\begin{equation} \label{transAinfty}
\mu \in A_\infty(\nu) \text{ and } \nu \in A_\infty(\sigma) \implies \mu \in A_\infty(\sigma).
\end{equation}
Moreover, we can adapt the definition to discuss \( A_\infty \)-absolute continuity of elliptic measures: we say that \( \omega_L^X \in A_\infty(\nu) \) if the implication \eqref{defAinfty} holds for all the measures \( \mu = \omega_L^X \) with \( X \in \Omega \setminus 2B \). 
\end{itemize}

\begin{proposition}[Lemma 3.2 in \cite{DFMpert}] \label{propGinfty}
Let $\Omega$ be an unbounded 1-sided Chord-Arc Domain, and let $L = -\diver A \nabla$ be a uniformly elliptic operator on $\Omega$ with adjoint $L^* = -\diver A^T \nabla$. Then there exist a function $G_{L^*}$ on $\Omega$ and a measure $\omega_L$ on $\partial \Omega$ such that:
\begin{enumerate}[(i)]
\item $G_{L^*}\in W^{1,2}_{loc}(\Omega)$ is a non-negative weak solution to $L^*u = -\diver A^T \nabla u = 0$ in $\Omega$;
\item $G_{L^*}\in C^0(\overline{\Omega})$ and $G_{L^*} \equiv 0$ on $\partial \Omega$;
\item For any $\varphi \in C_0^\infty(\mathbb{R}^n)$,
\[ \iint_{\Omega} A^T \nabla G_{L^*} \cdot \nabla \varphi \, dX = - \int_{\partial \Omega} \varphi(y) \, d\omega_L(y).\]
\end{enumerate}
The function $G_{L^*}$ and the measure $\omega_L$ are called the Green function with pole at infinity and the elliptic measure with pole at infinity, respectively. They are unique up to multiplication by a positive constant.
\end{proposition}

In the following, to simplify the notation, we write $\delta(X)$ for $\dist(X, \partial \Omega)$ and $B_X$ for the Whitney ball $B(X, \delta(X)/2)$.

\begin{definition} \label{defGL}
Let $\Omega$ be a bounded or unbounded 1-sided Chord-Arc Domain. We define the function $G_{L^*}$ and the measure $\omega_L$ as follows:
\begin{itemize}
\item[--] When $\Omega$ is unbounded, we take any $G_{L^*}$ and $\omega_L$ satisfying Proposition \ref{propGinfty};
\item[--] When $\Omega$ is bounded, we arbitrarily pick a point $X_0 \in \Omega$ such that $\delta(X_0) \geq \epsilon \diam(\Omega)$, where $\epsilon$ is the constant appearing in the definition of the Corkscrew point condition (see Subsection \ref{SsCAD}). We then set
\[ \omega_L := \omega^{X_0}_L \quad \text{and} \quad G_{L^*}(X) = \frac{1}{|B_{X_0}|} \iint_{B_{X_0}} G_{L^*}^Y(X) \, dY.\]
\end{itemize}
\end{definition}

The $A_\infty$-absolute continuity of the elliptic measure is related to the solvability of the Dirichlet problem in $L^p$, as shown in the following theorem.

\begin{theorem} \label{ThDPAinfty}
Let $\Omega$ be a 1-sided chord-arc domain with an Ahlfors regular measure $\sigma$, and let $L$ be a uniformly elliptic operator. Then the following are equivalent:
\begin{enumerate}[(i)]
\item $\omega_L^X \in A_\infty(\sigma)$;
\item $\omega_L \in A_\infty(\sigma)$;
\item There exists $p \in (1, \infty)$ such that the Dirichlet problem for the operator $L$ is solvable in $L^p(\partial \Omega, \sigma)$. That is, there exists $C > 0$ such that for any $g \in C^0_0(\partial \Omega)$, the function $u_g$ defined as in \eqref{defhm} satisfies the non-tangential bound
\[ \|N(u_g)\|_{L^p(\partial \Omega, \sigma)} \leq C \|g\|_{L^p(\partial \Omega, \sigma)},\]
where
\[ N(u_g)(x) := \sup_{\substack{X \in \Omega \\ |X - x| \leq 2\delta(X)}} |u_g(X)|.\]
\end{enumerate}
\end{theorem}

\begin{proof}The equivalence $(i) \iff (ii)$ follows from the comparison principle between the Green function and the elliptic measure, as demonstrated in Lemma~3.5 of \cite{DFMpert}.

The equivalence between $(i)$ and $(iii)$ in a stronger setting is established in Theorem~1.7.3 of \cite{KenigBook}, however this equivalence in 1-sided Chord-Arc Domains can be derived using a similar proof. Moreover, for further insights, one may look at Proposition~4.5 in \cite{HL}, which shows that, for the Laplacian on domains with Ahlfors regular boundaries, a weaker version of $(i)$ - specifically, the weak $A_\infty$ condition of the harmonic measure with respect to $\sigma$ - implies $(iii)$. Furthermore, Proposition~19 in \cite{Hof} demonstrates that $(iii)$ implies that $\omega^X_{-\Delta}$ is weak $A_\infty(\sigma)$. To complete the argument, it is important to note that in 1-sided Chord-Arc Domains, the equivalence between the weak $A_\infty$ and $A_\infty$ conditions of the elliptic measure is a direct consequence of the doubling property of the elliptic measure, see for instance Lemma~15.43 of \cite{DFMmixed}.
\end{proof}

\subsection{Carleson perturbations} 
Recall that \(\delta(X)\) denotes \(\dist(X, \partial \Omega)\), and \(B_X\) stands for \(B(X, \delta(X)/2)\). Observe that \(4B_X\) and \(8B_X\) are then ``Carleson regions,'' meaning that we can find \(x \in \partial \Omega\) and \(r > 0\) such that
\[ B(x, r) \subset 4B_X \subset 8B_X \subset B(x, 5r). \]

We first need a notion of a Carleson measure.

\begin{definition}
Let \(\Omega\) be a 1-sided Chord-Arc Domain.
\begin{itemize}
\item Let \(\sigma\) be an Ahlfors regular measure on \(\partial \Omega\). For a function \(f \in L^2_{loc}(\Omega)\), we define the norm
\[ \| f \|_{KCM(\sigma)} := \sup_{Y \in \Omega} \left( \frac{1}{\sigma(8B_Y \cap \partial \Omega)} \iint_{8B_Y \cap \Omega} |f(Z)|^2 \frac{dZ}{\delta(Z)} \right)^{\frac{1}{2}}.\]

\item Let \(L\) be a uniformly elliptic operator on \(\Omega\) and define \(G_{L^*}\) and \(\omega_L\) as in Definition \ref{defGL}. For a function \(f \in L^2_{loc}(\Omega)\), we define the norm
\[ \| f \|_{KCM(\omega_L)} := \sup_{Y \in \Omega} \left( \frac{1}{\omega_L(8B_Y \cap \partial \Omega)} \iint_{8B_Y \cap \Omega} |f(Z)|^2 G_{L^*}(Z) \frac{dZ}{\delta(Z)^2} \right)^{\frac{1}{2}}.\]
\end{itemize}
\end{definition}

\begin{proposition}[Theorem 1.29 in \cite{FP}] \label{prKCM}
Let \(\Omega\) be a 1-sided chord-arc domain, and let \(\mu, \nu\) be either an Ahlfors regular measure on \(\partial \Omega\), or an elliptic measure of a uniformly elliptic operator as given by Definition \ref{defGL}.

If \(\mu \in A_\infty(\nu)\), then there exists \(C > 0\) such that, for any \(f \in L^2_{loc}(\Omega)\),
\[ C^{-1} \| f \|_{KCM(\mu)} \leq \| f \|_{KCM(\nu)} \leq C \| f \|_{KCM(\mu)}.\]
In particular, \(\| f \|_{KCM(\mu)}\) is finite if and only if \(\| f \|_{KCM(\nu)}\) is finite.
\end{proposition}

The typical result on Carleson perturbations is of the following form (see \cite{FKP91, MPT, CHM, CHMT}):

\begin{theorem}[Theorem 1.3 in \cite{CHMT}] \label{ThFKP}
Let \(\Omega \subset \mathbb{R}^n\) be a 1-sided Chord-Arc Domain with Ahlfors regular measure \(\sigma\). Let \(L_i = -\diver A_i \nabla\), \(i \in \{0, 1\}\), be two uniformly elliptic operators and write \(\omega_i := \omega_{L_i}\) for the elliptic measure given by Definition \ref{defGL}. Define
\begin{equation}
\E_\infty(Y) := \sup_{Z \in B_Y} |A_1(Z) - A_0(Z)|.
\end{equation}
Assume finally that \(\|\E_\infty\|_{KCM(\sigma)} < +\infty\). Then \(\omega_0 \in A_\infty(\sigma)\) if and only if \(\omega_1 \in A_\infty(\sigma)\).
\end{theorem}


One might wonder if Theorem \ref{ThFKP} can be weakened by assuming a Carleson measure condition on an \(L^2\), or \(L^r\) average for finite \(r\), on the disagreement \(A_1 - A_0\) (instead of \(L^\infty\)). This involves replacing the quantity \(\E_\infty\) in Theorem \ref{ThFKP} with \(\E_2\) or \(\E_r\), where
\begin{equation} \label{defEr}
\E_r(Y) := \left(\frac{1}{|B_Y|} \iint_{Z \in B_Y} |A_1(Z) - A_0(Z)|^r \, dZ \right)^{\frac{1}{r}}.
\end{equation}

It has been shown in \cite[Proposition 2.22]{FKP91} (see also \cite{DSU}) that there exists an \(r \in [2, \infty)\), depending on the operator \(L_0\), such that \(\E_\infty\) can be replaced by \(\E_r\) in Theorem \ref{ThFKP}. Can we extend this to \(r = 2\)? While this may not hold in complete generality, it is possible under additional assumptions on the operator. Specifically, if the diffusion behavior induced by either \(L_0\) or \(L_1\) is uniformly ``nice" in the Whitney balls \(\{B_X\}_{X \in \Omega}\), then \(\E_2\) can be used instead of \(\E_\infty\).

Our main result is the following:

\begin{theorem} \label{MainTh}
Let \(\Omega \subset \mathbb{R}^n\) be a 1-sided CAD.  Let \(L_i = -\diver A_i \nabla\), \(i \in \{0, 1\}\), be two uniformly elliptic operators, and \(\omega_i := \omega_{L_i}\) be their elliptic measures in the sense of Definition \ref{defGL}. Assume finally that:
\begin{itemize}
\item[--] \(\|\E_2\|_{KCM(\omega_0)}\) is finite;
\item[--] either \(\delta \nabla A_0\) or \(\delta \nabla A_1\) belongs to \(L^\infty(\Omega)\).
\end{itemize}
Then $\omega_1 \in A_\infty(\omega_0)$.
\end{theorem}

\begin{corollary} \label{MainCor}
Let \(\Omega \subset \mathbb{R}^n\) be a 1-sided CAD with Ahlfors regular measure $\sigma$.  Let \(L_i = -\diver A_i \nabla\), \(i \in \{0, 1\}\), be two uniformly elliptic operators whose elliptic measures \(\omega_i := \omega_{L_i}\) are as in Definition \ref{defGL}. Assume finally that:
\begin{itemize}
\item[--] \(\|\E_2\|_{KCM(\sigma)}\) is finite;
\item[--] either $\delta \nabla A_0$ or $\delta \nabla A_1$ belongs to $L^\infty(\Omega)$.
\end{itemize}
Then $\displaystyle \omega_0 \in A_\infty(\sigma) \iff \omega_1 \in A_\infty(\sigma)$.
\end{corollary}

\begin{proof}[Proof of the Corollary]
The roles of $L_0$ and $L_1$ being symmetric, we just prove one implication. If \(\|\E_2\|_{KCM(\sigma)} < +\infty\) and \(\omega_0 \in A_\infty(\sigma)\), then \(\|\E_2\|_{KCM(\omega_0)} < +\infty\) by Proposition \ref{prKCM}. We can then apply Theorem \ref{MainTh} to conclude that $\omega_1 \in A_\infty(\omega_0)$, hence $\omega_1 \in A_\infty(\sigma)$ by transitivity of the $A_\infty$-absolute continuity.
\end{proof}

As an application of the \(L^2\) Carleson perturbation, we have the following result:

\begin{theorem} \label{ThHMMTZ}
Let \(n \geq 3\) and \(\Omega \subset \mathbb{R}^n\) be a 1-sided chord-arc domain. Let \(L = -\diver A \nabla\) be a uniformly elliptic operator and take the elliptic measure \(\omega := \omega_L\) as in Definition \ref{defGL}. Define
\[
\mathcal{O}_2(Y) = \inf_{C} \left(\fint_{B_Y} |A(Z) - C|^2 \, dZ\right)^{\frac{1}{2}},
\]
where the infimum is taken over the constant matrices. Assume that
\[
\|\mathcal{O}_2\|_{KCM(\sigma)} < +\infty.
\]
Then the following are equivalent:
\begin{enumerate}
\item \(\omega \in A_\infty(\sigma)\);
\item \(\partial \Omega\) is uniformly rectifiable\footnote{A quantitative version of rectifiability; see, for instance, \cite{DS1,DS2}.}.
\end{enumerate}
\end{theorem}

\begin{proof}[Proof of Theorem \ref{ThHMMTZ}, using \cite{HMMTZ}]
Take \(\theta \in C^\infty_0(\mathbb{R}^n)\) such that \(\supp \theta \subset B(0, \frac{1}{8})\), \(\iint_{\mathbb{R}^n} \theta(Z) \, dZ = 1\), and then construct \(\theta_Y(Z) := \delta(Y)^{-n} \theta\left(\frac{Z - Y}{\delta(Y)}\right)\) and
\[
B(Y) := \iint_{\mathbb{R}^n} A(Z) \theta_Y(Z) \, dZ.
\]
Observe that \(B\) is essentially the average of \(A\) over Whitney regions. We follow the proof of \cite[Corollary 2.3]{DPP} or \cite[Lemma 2.1]{FLM}, and we get that:
\begin{enumerate}[(a)]
\item \(\|\delta \nabla B\|_{L^\infty(\Omega)} < +\infty\) and \(\|\delta \nabla B\|_{KCM(\sigma)} < +\infty\),
\item \(\|A - B\|_{KCM(\sigma)} < +\infty\).
\end{enumerate}
The assertion (a) means that \(L_0 := -\diver B \nabla\) is a so-called DKP\footnote{For Dahlberg-Kenig-Pipher.} operator. So Theorem 1.6 in \cite{HMMTZ} gives that
\[
\omega_0 \in A_\infty(\sigma) \iff \partial \Omega \text{ is uniformly rectifiable.}
\]
Assertion (b) entails that \(L\) is an \(L^2\) Carleson perturbation of \(L_0\).  Indeed,  a simple application of Fubini's lemma yields that \(\|A - B\|_{KCM(\sigma)}\) is equivalent to \(\|\E_2\|_{KCM(\sigma)}\), where \(\E_2\) is the \(L^2\)-averaged disagreement between \(A\) and \(B\). Since \(\delta |\nabla B| \in L^\infty(\Omega)\), Corollary \ref{MainCor} gives
\[
\omega \in A_\infty(\sigma) \iff \omega_0 \in A_\infty(\sigma).
\]
The theorem follows.
\end{proof}

\subsection{Plan of the article}
In Section~2, we review the definition of 1-sided Chord Arc Domains and discuss several properties of $G_{L^*}$ and $\omega_L$, as introduced in Definition~\ref{defGL}. Section~3 focuses on proving reverse H\"older estimates for the gradient of solutions, which are essential for Section~4, where we present the core of the proof of Theorem~\ref{MainTh}.

\medskip

In the remainder of the article, the notation \( A \lesssim B \) means that there exists a constant \( C \), independent of the key parameters, such that \( A \leq CB \). Furthermore, we use \( A \approx B \) to denote that both \( A \lesssim B \) and \( B \lesssim A \) hold. Moreover, to lighten the notation, we often write $\iint_\Omega f$ for $\iint_\Omega f \, dX$, i.e. the Lebesgue measure will not always be mentioned in integrals.

\medskip

\noindent {\bf Acknowlegments:} The author would like to thank the anonymous referee for their suggestions, which significantly improved the readability of the article and identified a gap in the proof of Theorem~\ref{ThHMMTZ} in an earlier version of the manuscript.

\section{Preliminaries}
\subsection{The setting}  \label{SsCAD} 

Let \(\Omega \subset \mathbb{R}^n\) be an open set.

\begin{itemize}
\item We say that \(\partial \Omega\) is \((n-1)\)-{\bf Ahlfors regular} if there exist a constant \(C > 0\) and a measure \(\sigma\) supported on \(\partial \Omega\) such that
\begin{equation} \label{defAR}
C^{-1} r^{n-1} \leq \sigma(B(x, r)) \leq C r^{n-1} \quad \text{for } x \in \partial \Omega, \, r \in (0, \diam \Omega).
\end{equation}

\item We say that \(\Omega\) satisfies the {\bf Corkscrew point condition} if there exists \(\epsilon = \epsilon_{CP} > 0\) such that for any \(r \in (0, 2\diam \Omega)\) and any \(x \in \partial \Omega\), we can find \(X = X_{x, r}\) such that \(B(X, \epsilon r) \subset \Omega\). Note that if \(\Omega\) is bounded, by taking \(r = \diam \Omega\), we can find a `center' point \(X_0\) in \(\Omega\) such that \(\delta(X) \geq \epsilon_{CP} \diam \Omega\).

\item We say that \(\Omega\) satisfies the {\bf Harnack chain condition} if, for any \(\Lambda \geq 1\), we can find a constant \(C_\Lambda > 0\) such that for any pair \((X, Y) \in \Omega^2\) satisfying \(|X - Y| \leq \Lambda \min\{\delta(X), \delta(Y)\}\), there is a rectifiable curve \(\gamma = \gamma_{X, Y}: [0, \ell(\gamma)] \to \Omega\) parametrized by arclength such that:
\begin{enumerate}[(i)]
\item \(\gamma(0) = X\) and \(\gamma(\ell(\gamma)) = Y\).
\item The length \(\ell(\gamma)\) of \(\gamma\) is bounded by \(C_\Lambda |X - Y|\).
\item For each \(t \in (0, \ell(\gamma))\), \(B(\gamma(t), C_\Lambda^{-1} \min\{\delta(X), \delta(Y)\}) \subset \Omega\).
\end{enumerate}

\item If a domain \(\Omega\) satisfies both the Corkscrew point condition and the Harnack chain condition, we say that \(\Omega\) is {\bf 1-sided NTA}. 1-sided NTA domains are also known as uniform domains. Equivalent characterizations of such domains can be found in \cite{Vai88}.

\item A {\bf 1-sided Chord-Arc Domain} or {\bf 1 sided CAD} is simply a 1-sided NTA domain with a \((n-1)\)-Ahlfors regular boundary. In this scenario, an Ahlfors regular measure for \(\Omega\) is a measure \(\sigma\) satisfying \eqref{defAR}.
\end{itemize}

\subsection{Elliptic theory} 

\begin{theorem} \label{ThGandom}
Let $\Omega$ be a 1-sided CAD and $L=-\diver A \nabla$ be a uniformly elliptic operator. Let $G_{L^*}$ and $\omega_L$ be the Green function and elliptic measure from Definition \ref{defGL}. There exists $C>0$ such that for any $X\in \Omega$,
\[ C^{-1} \delta(X)^{n-2} G_{L^*}(X) \leq \omega_L(4B_X) \leq \omega_L(8B_X) \leq C \delta(X)^{n-2} G_{L^*}(X) \]
More specifically, the measure $\omega_L$ is doubling, meaning that there exists $C>0$ such that  for any $X\in \Omega$ and any $\Lambda \geq 4$, 
\[ \omega_L(2\Lambda B_X) \leq C\omega_L(2\Lambda B_X) .\]
\end{theorem}

\bp
In the bounded case, the theorem is just a special case of the comparison between Green function and elliptic measure. The proof originate from \cite{CFMS}, and the result in a setting more general than 1-sided CAD can be found as Lemmas 15.28 and 15.43 in \cite{DFMmixed}.

In the unbounded case, the proof stays true because the $G_{L^*}$ and $\omega_L$ are limits where $Y\to \infty$ of a rescaled version of $G_{L^*}^Y$ and $\omega_L^Y$. As such the equivalences that held for finite pole extend to the case where the pole is at infinity; see Lemmas 3.3 and 3.4 in \cite{DFMpert} for the result.
\ep

\section{Reverse H\"older inequality on the gradient of solutions.}

Reverse Hölder inequalities for the gradient of solutions to uniformly elliptic operators are classical. Some classical results are presented below:

\begin{proposition} \label{prRHs}
Let \( L = -\operatorname{div}(A \nabla) \) be a uniformly elliptic operator on a domain \( \Omega \). Then:
\begin{enumerate}[(i)]
\item There exist \( s > 2 \) and \( C > 0 \) such that for any ball \( B \subset \Omega \) and any solution \( u \) to \( Lu = 0 \) in \( B \), we have
\begin{equation} \label{defRHs}
\left( \frac{1}{\left| \frac{1}{2} B \right|} \iint_{\frac{1}{2} B} |\nabla u|^s\, dZ \right)^{\frac{1}{s}} 
\leq C \left( \frac{1}{|B|} \iint_{B} |\nabla u|^2\, dZ \right)^{\frac{1}{2}}.
\end{equation}
\item If, in addition, \( \delta \nabla A \in L^\infty(\Omega) \), then there exists \( C > 0 \) such that for any ball \( B \subset \Omega \) and any solution \( u \) to \( Lu = 0 \) in \( B \), we have
\begin{equation} \label{defRHinfty}
\sup_{\frac{1}{2} B} |\nabla u| 
\leq C \left( \frac{1}{|B|} \iint_{B} |\nabla u|^2\, dZ \right)^{\frac{1}{2}}.
\end{equation}
\end{enumerate}
The exponent \( s \) in (i) depends only on the dimension \( n \) and the ellipticity and boundedness constants of \( L \). In (i), the constant \( C \) in \eqref{defRHs} depends on the same parameters as \( s \), while in (ii), the constant \( C \) in \eqref{defRHinfty} also depends on \( \|\delta \nabla A\|_{L^\infty(\Omega)} \).
\end{proposition}

\begin{proof}
Case (i) is classical; see, for instance, Lemma 1.1.12 in \cite{KenigBook}. Case (ii) is also a straightforward consequence of standard results. Indeed, let \( r \) be the radius of \( B \), and denote by \( u_Y := \fint_{B(Y, r/4)} u \) the average of \( u \) over \( B(Y, r/4) \). Then, Lemma 3.1 in \cite{GW} yields
\[
|\nabla u(Y)| = |\nabla (u - u_Y)(Y)| 
\lesssim \frac{1}{r} \sup_{B(Y, r/4)} |u - u_Y| 
\lesssim \frac{1}{r} \fint_{B(Y, r/2)} |u - u_Y| 
\lesssim \left( \fint_{B(Y, r/2)} |\nabla u|^2 \right)^{\frac{1}{2}}.
\]
This follows from the Moser estimate and the Poincaré inequality. Taking the supremum over \( Y \in \frac{1}{2} B \), we obtain
\[
\sup_{\frac{1}{2} B} |\nabla u| 
\lesssim \sup_{Y \in \frac{1}{2} B} \left( \fint_{B(Y, r/2)} |\nabla u|^2 \right)^{\frac{1}{2}} 
\lesssim \left( \fint_{B} |\nabla u|^2 \right)^{\frac{1}{2}},
\]
as required.
\end{proof}

We will later require an improved version of the reverse Hölder estimate \eqref{defRHinfty} that incorporates cutoff functions. Moser-type estimates adapted to include cutoff functions have already been studied in Lemmas 3.1 and 3.3 of \cite{FMZ}, and we will follow their approach.

\begin{theorem} \label{ThRH}
Let $\Omega \subset \R^n$ be open and let $L=-\diver A \nabla$ be a uniformly elliptic operator on $\Omega$ that satisfies $\delta \nabla A \in L^\infty(\Omega)$. 

For any $s\in (2,\infty)$, there exists $k_s >0$ such that for any $k>k_s$, there exists $C_s>0$ such that for any ball $B$ satisfying $4B\subset \Omega$, any solution $u$ to $Lu = 0$ in $4B$,  and any cutoff function $\Psi$ satisfying $0 \leq \Psi \leq 1$ and $|\nabla \Psi| \leq 100/\delta$, we have
\[ \left(\frac1{|B|}\iint_B |\nabla u|^{s} \Psi^k \, dX\right)^\frac1q \leq C_{s,k} \left(\frac1{|2B|}\iint_{2B} |\nabla u|^{2} \Psi^{2k/s-k_s} \, dX\right)^\frac12.\]
\end{theorem}

\begin{proof}
\noindent\textbf{Step 0: A reminder on \( q \)-ellipticity.} The notion of \( q \)-ellipticity for \( q \in (1, \infty) \) was introduced in \cite{CD}; see also \cite{CM} for an earlier, related notion called \( q \)-dissipativity, and \cite{DP} for applications in the PDE setting. The precise definition of \( q \)-ellipticity is not needed here and we list only the properties relevant to our argument:

\begin{enumerate}
\item If \( A \) is a real-valued uniformly elliptic matrix, then $A$ is \( q \)-elliptic for all \( q \in (1, \infty) \).

\item \( A \) being \( q \)-elliptic implies that there exists \( \lambda_q > 0 \) such that for any non-negative \( \chi \in L^\infty(\Omega) \), and any real-valued function \( u \) with \( |u|^{q-2} |\nabla u|^2 \chi \in L^1(\Omega) \), we have
\begin{equation} \label{defqellip}
\iint_{\Omega} A \nabla u \cdot \nabla\left[|u|^{q-2} u\right] \chi \geq \lambda_q \iint_{\Omega} |u|^{q-2} |\nabla u|^2 \chi.
\end{equation}

\item Moreover, there exists a constant \( C_q > 0 \) such that for any function \( u \) satisfying \( |u|^{q/2} \in W^{1,2}_{\mathrm{loc}}(\Omega) \) and any point \( X \) where \( u(X) \neq 0 \), we have
\begin{equation} \label{prqellip}
C_q^{-1} |u(X)|^{q-2} |\nabla u(X)|^2 \leq \left| \nabla\left(|u|^{q/2 - 1} u\right)(X) \right|^2 \leq C_q |u(X)|^{q-2} |\nabla u(X)|^2.
\end{equation}
\end{enumerate}

Property (1) follows from \cite[Proposition 5.15]{CD}, while Properties (2) and (3) are given by \cite[Theorem 2.4 and Lemma 2.5]{DP}, respectively.

\medskip

\noindent\textbf{Step 1: The Caccioppoli inequality.} Let \( q > 2 \) and let \( u \) be a solution in a ball \( 3B \subset \Omega \). For any coordinate direction, let \( v_i = \partial_i u \). Since $\delta \nabla A \in L^\infty$, we have 
\begin{itemize}
\item[--] \( v_i \in L^\infty(2B) \cap W^{1,2}(2B) \), by \cite[Lemma 4.39]{HMT17};
\item[--] \( v_i \) is a weak solution to \( L v_i = -\operatorname{div}[(\partial_i A) \nabla u] \).
\end{itemize}

Let \( \eta_B \in C_0^\infty(2B) \) be a cutoff function such that \( \eta_B \equiv 1 \) on \( B \) and \( |\nabla \eta_B| \leq 2/r_B \), where \( r_B \) is the radius of \( B \). Define \( \Phi_B := \Psi \eta_B \). We claim that
\begin{equation} \label{claimCaccio}
I := \iint_\Omega |\nabla v_i|^2 |v_i|^{q-2} \Phi_B^k \lesssim J := \frac{1}{r^2_B} \iint_\Omega |\nabla u|^q \Phi_B^{k-2},
\end{equation}
with constants depending on \( q \) and \( k \), but independent of \( B \) and \( u \).

Since \( v_i \in L^\infty(2B) \cap W^{1,2}(2B) \), the integrals in \eqref{claimCaccio} are finite. For simplicity, we now write \( v := v_i \), \( \Phi := \Phi_B \), and \( r := r_B \). Using the \( q \)-ellipticity inequality \eqref{defqellip}, we estimate
\begin{align*}
I &\lesssim \iint_\Omega A \nabla v \cdot \nabla \left[|v|^{q-2} v\right] \Phi^k \\
&= \iint_\Omega A \nabla v \cdot \nabla \left[|v|^{q-2} v \Phi^k \right] - \iint_\Omega A \nabla v \cdot \nabla \Phi^k |v|^{q-2} v \\
&=: I_1 + I_2.
\end{align*}
Because \( v \) is a weak solution to \( Lv = -\operatorname{div}[(\partial_i A) \nabla u] \), we get
\begin{align*}
|I_1| &= \left| \iint_\Omega (\partial_i A) \nabla u \cdot \nabla \left[|v|^{q-2} v \Phi^k\right] \right| \\
&\lesssim \frac{\|\delta \nabla A\|_{\infty}}{r} \iint_\Omega |\nabla u| \left( \Phi^k |\nabla [|v|^{q-2} v]| + |v|^{q-1} |\nabla \Phi^k| \right) \\
&\lesssim \frac{1}{r} \iint_\Omega |\nabla u| |v|^{q-2} \Phi^{k-1} \left( \Phi |\nabla v| + |v| |\nabla \Phi| \right) \\
&\lesssim \iint_\Omega |v|^{q-2} \Phi^{k-1} \left( \Phi \frac{|\nabla u|^2}{\epsilon r^2} + \epsilon \Phi |\nabla v|^2 + |v| |\nabla u| \frac{|\nabla \Phi|}{r} \right),
\end{align*}
where we used \( |\nabla A| \lesssim \delta^{-1} \lesssim 1/r \) in \( 2B \) and Young’s inequality in the last line. Since \( |v| \leq |\nabla u| \), \( \Phi \leq 1 \), and \( |\nabla \Phi| \lesssim 1/r \), we conclude that
\[
I_1 \leq C \left( \epsilon I + \frac{1}{\epsilon} J \right), \quad \text{for any } \epsilon \in (0,1).
\]
Similarly,
\[
|I_2| \lesssim \iint_\Omega |\nabla v| \Phi^{k-1} |\nabla \Phi| |v|^{q-1} \lesssim \epsilon I + \frac{1}{\epsilon} J,
\]
since \( |\nabla \Phi| \lesssim 1/r \). Altogether, we obtain
\[
I \lesssim \epsilon I + \frac{1}{\epsilon} J.
\]
Choosing \( \epsilon \) small enough and using that \( I < \infty \), we deduce \( I \lesssim J \), establishing \eqref{claimCaccio}.

\medskip

\noindent\textbf{Step 2: The reverse Hölder inequality.} We now combine the Caccioppoli inequality with a Sobolev-Poincaré inequality to derive a reverse Hölder estimate.

Let \( q \leq p \leq \frac{n}{n-1} q \). For \( k \geq 2n/(n-1) \geq 2p/q \), define \( w := |v|^{q/2} \Phi^{kq/2p} \in W^{1,2}(2B) \), compactly supported in \( 2B \). Then, by the $L^{2p/q}$-$L^2$ Sobolev-Poincar\'e inequality, we have
\begin{align*}
\left( \frac1{|2B|}\iint_{2B} |v|^p \Phi^k \right)^{1/p} 
&= \left( \frac1{|2B|}\iint_{2B} w^{2p/q} \right)^{1/p} 
\lesssim \left( \frac{r^2}{|2B|} \iint_{2B} |\nabla w|^2 \right)^{1/q} \\
&\lesssim \left( \frac1{|2B|}\iint_{2B} |\nabla v|^2 |v|^{q-2} \Phi^{kq/p} \right)^{1/q} 
+ \left( \frac1{|2B|}\iint_{2B} |v|^q \Phi^{kq/p - 2} (r |\nabla \Phi|)^2 \right)^{1/q}.
\end{align*}
Using \eqref{claimCaccio} and \( r|\nabla \Phi| \lesssim 1 \), we conclude:
\begin{equation} \label{wRHb}
\left( \frac1{|2B|}\iint_{2B} |v|^p \Phi^k \right)^{1/p} \lesssim \left( \frac1{|2B|}\iint_{2B} |\nabla u|^q \Phi^{kq/p - 2} \right)^{1/p}.
\end{equation}
Since \( |\nabla u|^q \approx \sum_{i=1}^n |v_i|^q \), we obtain
\begin{equation} \label{wRH}
\left(\frac1{|2B|}\iint_{2B} |\nabla u|^p \Phi^k \right)^{1/p} \lesssim \left( \frac1{|2B|}\iint_{2B} |\nabla u|^q \Phi^{kq/p - 2} \right)^{1/p}.
\end{equation}

Fix \( s > 2 \). Then there exists \( j_s \in \mathbb{N} \) such that \( s/2 \in \left( \left( \frac{n}{n-1} \right)^{j-1}, \left( \frac{n}{n-1} \right)^j \right] \). Let \( k_s := 2j \). Iterating \eqref{wRH}, we deduce:
\begin{equation} \label{sRH}
\left( \frac1{|2B|}\iint_{2B}  |\nabla u|^s \Phi^k \right)^{1/s} \lesssim \left(\frac1{|2B|}\iint_{2B} |\nabla u|^2 \Phi^{2k/s - k_s} \right)^{1/2}.
\end{equation}
Removing the cutoff \( \eta_B \), we conclude:
\[
\left( \frac1{|B|}\iint_{B} |\nabla u|^s \Psi^k \right)^{1/s} \lesssim \left(\frac1{|2B|}\iint_{2B} |\nabla u|^2 \Psi^{2k/s - k_s} \right)^{1/2},
\]
which completes the proof.
\end{proof}

\section{Proof of Theorem \ref{MainTh}}

The absolute continuity between the elliptic measures will be a consequence of the following theorem, whose idea originates in \cite{KKPT16}, that was adapted in various settings in \cite{DFMDahlberg,CHMT,FP,AHMT23}.

\begin{theorem}[{\cite[Theorem 1.3]{AHMT23}, \cite[Theorem 1.22]{FP}}] \label{ThCME=>Ainfty}
Let $\Omega \subset \R^n$ be a $1$-sided CAD. For $i\in \{0,1\}$, let $L_i=-\diver A_i \nabla$ be two uniformly elliptic operators, and write $G_i = G_{L^*_i}$ and $\omega_i := \omega_{L_i}$ for the Green function and the elliptic measure from Definition \ref{defGL}. Assume that there exists $C>0$ such that any bounded weak solution $u$ to $L_1u=0$ in $\Omega$ and any $X\in \Omega$, we have
\begin{equation} \label{CME}
\iint_{4B_X \cap \Omega} |\nabla u|^2 \, G_0\, dZ \leq C \|u\|_\infty^2 \omega_{0}(8B_X \cap \partial \Omega),
\end{equation}
Then $\omega_1 \in A_\infty(\omega_0)$.
\end{theorem}

To prove Theorem~\ref{MainTh}, it suffices to verify the assumptions of Theorem~\ref{ThCME=>Ainfty}, specifically \eqref{CME}. The strategy employed is classical, and we follow closely the proof of Lemma~5.1 in \cite{FP}.

\medskip

\textbf{Step 1: Setting the bound.} Let \( X \in \Omega \) be given, and let \( u \) be a bounded weak solution to \( L_1u = 0 \) in \( \Omega \). We construct a smooth cutoff function \( \phi \) such that \( \phi \equiv 1 \) on \( 4B_X \), \( \phi \equiv 0 \) on \( \mathbb{R}^n \setminus 7B_X \), and \( |\nabla \phi| \lesssim \delta(X) \). For a \( k \geq 2 \) (independent of \( X \) and \( u \)) to be determined later, we aim to estimate
\[
\iint_{\Omega} |\nabla u|^2 \phi^k \, G_0,
\]
but since we do not know \emph{a priori} that such a quantity is finite, we fix \( i \in \mathbb{N} \) and use an additional cutoff function \( \psi_i \) such that \( \psi_i(Z) = 1 \) when \( \delta(Z) \geq 2^{-i} \), \( \psi_i(Z) = 0 \) when \( \delta(Z) \leq 2^{-i-1} \), and \( |\nabla \psi_i| \leq 2^{3-i} \). We seek to prove that
\begin{equation} \label{claim}
\iint_{\Omega} |\nabla u|^2 \phi^k \psi_i^k \, G_0 \leq C\|u\|_\infty^2 \omega_0(8B_X \cap \partial \Omega),
\end{equation}
with a constant \( C \) independent of \( X \), \( Y \), \( u \), and \( k \). Once \eqref{claim} is established, we take \( k \to 0 \) to obtain the bound
\[
\iint_{4B_X} |\nabla u|^2 G_0 \leq \iint_{\Omega} |\nabla u|^2 \phi^2 G_0 \leq C\|u\|_\infty^2 \omega_0(8B_X \cap \partial \Omega),
\]
required to invoke Theorem~\ref{ThCME=>Ainfty}.

\medskip

\textbf{Step 2: Carleson estimate on $\phi\psi_i$.} For the remainder of the proof, we denote $\varphi = \varphi_i$ for $\phi\psi_i$. We aim to show that there exists a constant $C > 0$, independent of $i$, such that
\begin{equation} \label{claimphi}
\iint_{\Omega} |\nabla \varphi|^2 G_0 + \iint_{\Omega} |\nabla \varphi| |\nabla G_0| \leq C \omega_0(8B_X \cap \partial \Omega).
\end{equation}
The comparison between the elliptic measure and the Green function (Theorem~\ref{ThGandom}) yields $\omega_0(8B_Z \cap \partial \Omega) \approx \delta(Z)^{2-n} G_0(Z)$. Consequently, we deduce that
\begin{multline*}
\iint_{\Omega} |\nabla \varphi|^2 G_0 \lesssim \iint_{8B_X} |\nabla \varphi(Z)|^2 \delta(Z)^{2-n} \omega_0(8B_Z \cap \partial \Omega) \, dZ \\
= \iint_{8B_X} |\nabla \varphi(Z)|^2 \delta(Z)^{2-n} \left( \int_{8B_Z \cap \partial \Omega} d\omega_0(y)\right) \, dZ \\
= \iint_{y \in H_X} \left( \underbrace{\int_{|Z-y| \leq 8 \delta(Z)} \delta(Z)^{2-n} |\nabla \varphi|^2 \, dZ}_{=:f(y)} \right) d\omega_0(y),
\end{multline*}
by Fubini's theorem, with $H_X := \bigcup_{Z \in 8B_X} 8B_Z$. Note first that $H_X \subset 80B_X$. We then aim to prove that $f \lesssim 1$. To this end, we define $\gamma^*_j(y) := \{Z \in \Omega, |Z-y| \leq 8 \delta(Z), 2^{-j-1} \leq \delta(Z) \leq 2^{-j}\}$ and observe that
\[ f(y) = \sum_{j \in \mathbb{Z}} \iint_{\gamma^*_j(y) \cap 8B_X} \delta^{2-n} |\nabla \varphi|^2 \leq \sum_{j \in \mathbb{Z}} 2^{(n-2)(j+1)} \iint_{\gamma^*_j(y) \cap 8B_X} |\nabla \varphi|^2. \]
Noticing that
\[ |\nabla \varphi(Z)| \lesssim \frac{1}{\delta(X)} + 2^i \mathbf{1}_{\{2^{-i-1} \leq \delta(Z) \leq 2^{-i}\}}, \]
and that $|\gamma^*_j(y)| \approx 2^{-jn}$, we deduce
\[ f(y) \lesssim 2^{2i} 2^{(n-2)i} |\gamma^*_k(y)| + \frac{1}{\delta(X)^2} \sum_{j \geq -\ln(2\delta(X))} 2^{(n-2)j} |\gamma^*_j(y)| \lesssim 1. \]
We conclude that
\[ \iint_{\Omega} |\nabla \varphi|^2 G_0 \, dZ \lesssim \omega_0(80B_X \cap \partial \Omega) \lesssim \omega_0(8B_X \cap \partial \Omega) \]
by the doubling property of $\omega_0$ (Theorem~\ref{ThGandom}).

It remains to prove the second part of the claim \eqref{claimphi}. This part is similar: thanks to the Cauchy-Schwarz inequality, the Caccioppoli inequality, the Harnack inequality, and the comparison between the Green function and the elliptic measure, we have, for any $Z \in \Omega$,
\[ \left( \iint_{\frac{1}{2} B_Z} |\nabla G_0|^2 dA \right)^{\frac{1}{2}} \lesssim \delta(Z)^{\frac{n}{2}-2} \left( \iint_{B_Z} |G_0|^2 dA \right)^{\frac{1}{2}} \lesssim \delta(Z)^{n-2} G(Z) \approx \omega^X_0(8B_Z \cap \partial \Omega). \]

Then, Fubini's theorem and the Cauchy-Schwarz inequality give that
\begin{multline*}
\iint_{\Omega} |\nabla \varphi| |\nabla G_0| = \iint_{\Omega} \left( \frac{1}{|\frac{1}{2} B_Z|} \iint_{\frac{1}{2} B_Z} |\nabla \varphi| |\nabla G_0| \right) \, dZ \\
\lesssim \iint_{\Omega} \left( \frac{1}{|\frac{1}{2} B_Z|} \iint_{\frac{1}{2} B_Z} |\delta \nabla \varphi|^2 \right)^{\frac{1}{2}} \left( \iint_{\frac{1}{2} B_Z} |\nabla G_0|^2 \right)^{\frac{1}{2}} \, dZ \\
\lesssim \iint_{\Omega} \delta(Z)^{1-n} \left( \frac{1}{|\frac{1}{2} B_Z|} \iint_{\frac{1}{2} B_Z} |\nabla \varphi|^2 \right)^{\frac{1}{2}} \omega_0(8B_Z \cap \partial \Omega) \, dZ.
\end{multline*}
The rest of the argument follows the same lines as those used to bound $\iint_{\Omega} |\nabla \varphi|^2 G_0$.

\medskip

\noindent {\bf Step 3: Integration by part.}  Define 
\[I := \iint_\Omega |\nabla u|^2 \varphi^k G_0  \quad \text{ and } \quad J:= \|u\|_\infty^2 \omega_0(8B_X \cap \partial \Omega).\]
We want to prove that $I\lesssim J + (IJ)^{1/2}$, which implies the claim \eqref{claim} since $I$ is finite. By the ellipticity of $A_1$, we first have
\[ I \lesssim  \iint_{\Omega} A_1 \nabla u \cdot \nabla u\, \varphi^k \, G_0 \]
We want to use the fact that $u$ is a solution to $Lu_1 = 0$. Since $u\varphi G_0$ lies in $W^{1,2}(\Omega)$ and is compactly supported in $8B_Y \cap \Omega$, it is a valid test function, so 
\[ \iint_{\Omega} A_1 \nabla u \cdot \nabla [u\varphi^k G_0]  = 0.\]
We deduce that
\[I \lesssim \iint_{\Omega} A_1 \nabla u \cdot \nabla u\, \varphi^k \, G_0\\
 = - k\iint_{\Omega} A_1 \nabla u \cdot \nabla \varphi\, u \varphi^{k-1} \, G_0 - \iint_{\Omega} A_1 \nabla u \cdot \nabla G_0\, u \varphi^k =: I_1+ I_2. \]
The term $I_1$ is bounded using the Cauchy-Schwarz inequality (together with the facts that $A_1 \in L^\infty$ and $\varphi \leq 1$)
\[ |I_1| \lesssim I^\frac12 \|u\|_\infty \left(  \iint_\Omega |\nabla \varphi|^2 \varphi^{k-2} G_0 \right)^\frac12 \lesssim (IJ)^\frac12\]
by \eqref{claimphi}. As for the term $I_2$, we want to use the fact that $G_0$ is a weak solution for the operator $L_0^* := -\diver A_0^T \nabla$. So we write
\begin{multline*} 
I_2 = - k\iint_\Omega A_0 \nabla \varphi \cdot \nabla G_0 \, u^2 \varphi^{k-1} -\frac 12 \iint_{\Omega} A_0 \nabla[u^2\varphi^k] \cdot \nabla G_0  \\  + \iint_{\Omega} (A_0-A_1) \nabla u \cdot \nabla G_0 \, u \varphi^k  =: I_3 + I_4 + I_5.
\end{multline*}
To estimate \( I_3 \), we simply invoke the second part of \eqref{claimphi}. Since \( A_0 \) is bounded, we obtain
\[
|I_3| \lesssim \|u\|_\infty^2 \left( \iint_\Omega |\nabla \varphi|\, |\nabla G_0| \right) \lesssim J.
\]
To handle \( I_4 \), we use the fact that \( G_0 \) is a solution. When the domain \( \Omega \) is unbounded, \( G_0 \) is a weak solution to \( L_0^* u = 0 \), and thus \( I_4 = 0 \). In the case where \( \Omega \) is bounded, let \( X_0 \) be the center used to define \( G_0 \). Then,
\[
I_4 = \frac{1}{|B_{X_0}|} \iint_{B_{X_0}} u^2 \varphi^2.
\]
We distinguish two cases:
\begin{itemize}
\item[--] If \( 8B_X \cap B_{X_0} = \emptyset \), then \( I_4 = 0 \).
\item[--] If \( 8B_X \cap B_{X_0} \neq \emptyset \), then \( |X - X_0| \lesssim \delta(X) \approx \delta(X_0) \), which yields, by the doubling property of \( \omega_0 \) (see Theorem~\ref{ThGandom}),
\[
\omega_0(8B_X \cap \partial \Omega) \approx \omega_0(8B_{X_0} \cap \partial \Omega) \approx \omega_0(\partial \Omega) = 1,
\]
since \( \omega_0 \) is a probability measure. To justify the second estimate above, observe that if \( \delta(X_0) \geq \epsilon \operatorname{diam}(\Omega) \), then \( 2\epsilon^{-1} B_{X_0} \supset \partial \Omega \), and so the doubling property of \( \omega_0 \) implies \( \omega_0(8B_{X_0}) \approx \omega_0(\partial \Omega) \).
\end{itemize}
In all cases, we conclude
\[
|I_4| \lesssim \|u\|_\infty^2 \lesssim \|u\|_\infty^2 \omega_0(8B_X \cap \partial \Omega) = J.
\]

We conclude with the estimate for \( I_5 \), distinguishing again two cases:

\smallskip

\noindent {\textbullet} \textbf{Case 1:} \( \delta |\nabla A_0| \in L^\infty \). By Proposition~\ref{prRHs}, we have \( |\nabla G_0| \lesssim G_0 / \delta \). It follows that
\begin{multline*}
|I_5| \lesssim \|u\|_\infty \iint_{\Omega} |A_0 - A_1|\, |\nabla u|\, \frac{G_0}{\delta} \varphi^k 
\leq \|u\|_\infty\, I^{1/2} \left( \iint_\Omega |A_0 - A_1|\, \frac{G_0}{\delta^2} \varphi^k \right)^{1/2} \\
\lesssim \|u\|_\infty\, I^{1/2}\, \omega_0(8B_X \cap \partial \Omega)\, \|A_0 - A_1\|_{CM(\omega_0)},
\end{multline*}
by Hölder's inequality and the fact that \( \operatorname{supp} \varphi \subset 8B_X \). Since Fubini's theorem implies \( \|A_0 - A_1\|_{CM(\omega_0)} \approx \|\mathcal{E}_2\|_{CM(\omega_0)} < \infty \), we deduce
\[
|I_5| \lesssim (I J)^{1/2},
\]
as desired.

\smallskip

\noindent {\textbullet} \textbf{Case 2:} \( \delta |\nabla A_1| \in L^\infty \).  Applying Fubini's theorem, we write
\[
|I_5| \lesssim \iint_{\Omega} \left( \frac{1}{|\frac{1}{2}B_Z|} \iint_{\frac{1}{2}B_Z} |A_0 - A_1|\, |\nabla u|\, |\nabla G_0|\, u\, \varphi^k \right)\, dZ.
\]
We aim to apply Hölder's inequality to the inner integral over \( \frac{1}{2}B_Z \). Proposition~\ref{prRHs} gives the existence of \( q > 2 \) such that \( \nabla G_0 \in L^q_{\text{loc}} \). Choosing \( s > 2 \) so that \( \frac{1}{q} + \frac{1}{s} = \frac{1}{2} \), Hölder's inequality yields
\[
|I_5| \lesssim \|u\|_\infty \iint_{\Omega} \mathcal{E}_2(Z) 
\left( \frac{1}{|\frac{1}{2}B_Z|} \iint_{\frac{1}{2}B_Z} |\nabla u|^s \varphi^{ks} \right)^{1/s}
\left( \frac{1}{|\frac{1}{2}B_Z|} \iint_{\frac{1}{2}B_Z} |\nabla G_0|^q \right)^{1/q} \, dZ.
\]
Choosing \( k = k_s \), where \( k_s \) is the exponent from Theorem~\ref{ThRH}, we apply the reverse Hölder inequality to get
\[
\left( \frac{1}{|\frac{1}{2}B_Z|} \iint_{\frac{1}{2}B_Z} |\nabla u|^s \varphi^{ks} \right)^{1/s}
\lesssim \left( \frac{1}{|B_Z|} \iint_{B_Z} |\nabla u|^2 \varphi^{2k - k_s} \right)^{1/2}
\leq \left( \frac{1}{|B_Z|} \iint_{B_Z} |\nabla u|^2 \varphi^k \right)^{1/2}.
\]
Moreover, using Proposition~\ref{prRHs}, the Caccioppoli inequality, and the Harnack inequality, we get
\[
\left( \frac{1}{|\frac{1}{2}B_Z|} \iint_{\frac{1}{2}B_Z} |\nabla G_0|^q \right)^{1/q}
\lesssim \left( \frac{1}{|B_Z|} \iint_{B_Z} |\nabla G_0|^2 \right)^{1/2}
\lesssim \frac{1}{\delta(Z)} \inf_{B_Z} G_0.
\]
Putting everything together,
\begin{multline*}
|I_5| \lesssim \|u\|_\infty \iint_{\Omega} \mathcal{E}_2(Z) \frac{G_0^{1/2}(Z)}{\delta(Z)}
\left( \frac{1}{|B_Z|} \iint_{B_Z} |\nabla u|^2 \varphi^k G_0 \right)^{1/2} \, dZ \\
\leq \|u\|_\infty 
\left( \underbrace{\iint_{50B_X} |\mathcal{E}_2|^2 \frac{G_0}{\delta^2}}_{\lesssim \omega_0(50B_X)} \right)^{1/2}
\left( \underbrace{\iint_\Omega \left( \frac{1}{|B_Z|} \iint_{B_Z} |\nabla u|^2 \varphi^k G_0 \right) \, dZ}_{\lesssim I \text{ by Fubini}} \right)^{1/2}
\lesssim (I J)^{1/2}.
\end{multline*}
For the second bound above, note that if \( B_Z \cap \operatorname{supp} \varphi \neq \emptyset \), then \( Z \in 50B_X \). The final bound follows from the finiteness of \( \|\mathcal{E}_2\|_{KCM(\omega_0)} \) and the doubling property of \( \omega_0 \). This completes the proof.

\bibliographystyle{amsalpha}

\end{document}